\documentclass[11.5pt]{article}
\usepackage{mathrsfs}

\usepackage{amsmath}
\usepackage{stmaryrd}
\usepackage{txfonts}
\usepackage{amsfonts}
\usepackage{cite}

\newcommand{\ba}{\noindent $\begin{array}}
\newcommand{\ea}{\end{array}$}
\newcommand{\be}{\begin{equation}}
\newcommand{\ee}{\end{equation}}
\newcommand{\bd}{\begin{displaymath}}
\newcommand{\ed}{\end{displaymath}}
\newcommand{\beq}{\begin{eqnarray*}}
\newcommand{\eeq}{\end{eqnarray*}}
\newcommand{\beqn}{\begin{eqnarray}}
\newcommand{\eeqn}{\end{eqnarray}}

\makeatletter

\newcommand{\Rmnum}[1]{\expandafter\@slowromancap\romannumeral #1@}
\makeatother

\usepackage{latexsym,lscape}

\newtheorem{theorem}{Theorem}[section]
\newtheorem{proposition}{Proposition}[section]
\newtheorem{lemma}{Lemma}[section]
\newtheorem{corollary}{Corollary}[section]
\newtheorem{definition}{Definition}[section]
\newtheorem{remark}{Remark}[section]

\newfont{\Bb}{msbm10 scaled\magstep1}

\def\sqr#1#2{{\vcenter{\hrule height .#2pt
      \hbox{\vrule width .#2pt height#1pt \kern#1pt\vrule width.#2pt}
                       \hrule height.#2pt}}}

\title{Optimization with Least Constraint Violation
}

\author{Yu-Hong Dai\footnote{LSEC, ICMSEC, AMSS, Chinese Academy of Sciences, Beijing 100190, China.  {\sl Email}: dyh@lsec.cc.ac.cn.
This author was supported by the Natural Science Foundation of China (Nos. 11991020, 12021001, 11631013, 11971372 and 11991021) and
Beijing Academy of Artificial Intelligence (BAAI).}
 \footnote{School of Mathematical Sciences, University of Chinese Academy of Sciences, Beijing 100049, China.}\quad and \quad Liwei Zhang \footnote{Corresponding author. School of Mathematical Sciences, Dalian University of Technology, Dalian 116024, China.\ {\sl Email}: lwzhang@dlut.edu.cn. This author was supported by the Natural Science Foundation of China (Nos. 11971089 and 11731013).}}

\begin{document}

\maketitle
\begin{abstract} Study about theory and algorithms for constrained optimization usually assumes that the feasible region of the optimization problem is nonempty. However, there are many important practical optimization problems whose feasible regions are not known to be nonempty or not, and  optimizers of the objective function with the least constraint violation prefer to be found.
A natural way for dealing with these problems is to extend the constrained optimization problem  as the one optimizing the objective function over the set of points with the least constraint violation.
  Firstly, the minimization problem with least constraint violation is proved to be an Lipschitz equality constrained optimization problem when the  original problem
is a convex  optimization problem with possible inconsistent conic constraints, and it can be reformulated as an MPEC problem. Secondly, for nonlinear programming problems with possible inconsistent constraints, various types of stationary points are presented for the MPCC problem which is equivalent to  the minimization problem with least constraint violation, and an elegant necessary optimality condition, named as L-stationary condition, is established from the classical optimality theory of Lipschitz continuous optimization.
Finally, the smoothing  Fischer-Burmeister function method for nonlinear programming case is constructed for  solving the  problem minimizing the objective function with the least constraint violation. It is demonstrated that, when the positive smoothing parameter approaches to zero, any point in the outer limit of
the KKT-point mapping is an L-stationary point of the equivalent MPCC problem.

\vskip 12 true pt \noindent \textbf{Key words}:  least constraint violation, L-stationary point, MPEC, MPCC, smoothing function.

\end{abstract}


\section{Introduction}
\setcounter{equation}{0}
For studying an ordinary nonlinear optimization problem, a basic assumption  is that feasible region of the optimization problem is nonempty. Many important theoretical issues are well studied for an optimization problem under this assumption.
For example, optimality theory and sensitivity analysis are two main theoretical topics. Optimality theory consists of necessary optimality conditions and sufficient optimality conditions.  Sensitivity analysis studies continuity properties of the optimal value and the solution mapping when the optimization is perturbed. For nonlinear programming, for a local minimizer, the first-order necessary optimality conditions and the second-order optimality conditions can be developed under certain constraint qualifications, and the second-order sufficient optimality conditions imply the second-order growth condition, see for instance the famous textbook \cite{NWright99}. For nonlinear programming, a series of stability results were obtained by Robinson, see \cite{Rob80},\cite{Rob81} and \cite{Rob82}.
 Bonnans and Shapiro \cite{BS00} established the optimality theory and the stability theory for general optimization problems, including problems whose decision variables are  infinite dimensional, nonlinear semidefinite programming problems and other conic optimization problems.

 However, when the feasible set is empty or the constraints are inconsistent, infeasibility detection is an important issue for algorithmic design. Many numerical algorithms have been proposed to find infeasible stationary points; namely, stationary points for minimizing certain infeasibility measure.  Byrd, Curtis and Nocedal \cite{Byrd2010} presented a set of conditions to guarantee the superlinear convergence of their SQP algorithm to an infeasible stationary
point.  Burke, Curtis and Wang \cite{Burke2014} considered the general program
with equality and inequality constraints, and proved that their SQP method has
strong global convergence and rapid convergence to the KKT point, and has
superlinear/quadratic convergence to an infeasible stationary point. Recently, Dai, Liu and Sun \cite{DLSun2020}
proposed a  primal-dual interior-point method, which
can be superlinearly or quadratically convergent to the Karush-Kuhn-Tucker point if the original problem is feasible, and can be superlinearly or quadratically convergent to the infeasible stationary point when the problem
is infeasible.

 These algorithms can find a stationary point of the infeasibility measure, which have nothing to do with the objective function of the problem. In practice, there are many important problems that we need to find  minimizers of the objective function over the points with the least constraint violation. A natural way to deal with such problems is to  extend the constrained optimization problem as the one that optimizes the objective function over the set of points with least constraint violation. When the the feasible region is nonempty, the set of points with least constraint violation coincides with the feasible region of the constrained optimization problem and hence the extended constrained optimization problem coincides with the original problem.

Now we give a formal definition of infeasibility measure of the optimization problem. Suppose that the problem is of the following form
\begin{equation}\label{p:ap}
\begin{array}{ll}
\min & f(x)\\[4pt]
{\rm s.t.} & c(x) \in K,
\end{array}
\end{equation}
where $f: {\cal X} \rightarrow \Re$, $c:{\cal X} \rightarrow {\cal Y}$, $K\subset {\cal Y}$, ${\cal X}$  and
${\cal Y}$ are finite dimensional Hilbert spaces.
\begin{definition}\label{def:olv}
A function $\theta:{\cal X}\rightarrow \Re$ is said to be an infeasibility measure of constraint $c(x) \in K$ if there exists an increasing continuous function $\varrho:\Re_+ \rightarrow \Re_+$ with $\varrho (0)=0$ such that
$$
\theta(x)=\varrho ({\rm dist }\,(c(x),K)),
$$
where
$$
{\rm dist }\,(c(x),K)=\inf\{\|y-c(x)\|:y\in K\}
$$
is the distance from $c(x)$ to $K$ under the norm $\|\cdot\|$ in ${\cal Y}$.
\end{definition}

Under the infeasibility measure defined above, we introduce the mathematical model of minimizing the objective $f(x)$ over the set of points with least infeasibility measure.
\begin{definition}\label{def:plc}
For  an infeasibility measure $\theta (x)$ of the constraint $c(x) \in K$, the mathematical model of minimizing the objective $f(x)$ over the set of points with least constraint violation associated with $\theta$, is defined by
\begin{equation}\label{def:mlv}
\left
\{
\begin{array}{cl}
\min & f(x)\\[5pt]
{\rm s.t.} & x \in {\rm Argmin}_z\, \theta (z).
\end{array}
\right.
\end{equation}
\end{definition}
Obviously, if the feasible region $c^{-1}(K)$ is nonempty, then $\min_z\, \theta (z)=0$, ${\rm Argmin}_z\, \theta (z)=c^{-1} (K)$, and Problem (\ref{def:mlv}) is just the original problem (\ref{p:ap}). Thus Problem (\ref{def:mlv}) can be regarded as an extension of  the original problem (\ref{p:ap}). As the minimization of the constraint violation is considered absolutely prior to the optimization of the objective function, we call the optimum (optimizer) of Problem (\ref{def:mlv}) as the {\sl optimum (optimizer) of Problem (\ref{p:ap}) with least constraint violation}, or simply {\sl constrained optimum (optimizer) of Problem (\ref{p:ap})}.

We do not identify the notion of ${\rm Argmin}_z\, \theta (z)$ in Definition \ref{def:mlv}. If $\theta (z)$ is convex, it is obvious that the solution set is just the set of global minimizers. However, if $\theta (z)$ is non-convex (this happens when  $\{x: c(x) \in K\}$ do not represent a convex set if feasible), ${\rm Argmin}_z\, \theta (z)$ may be understood as a set of local minimizers or even the set of stationary points. Infeasibility detection is well known  difficult in the nonconvex optimization case.
Indeed, for a nonconvex problem, infeasibility detection has many of the difficulties
inhered in global optimization since even if an algorithm identifies an infeasible point
where constraint violations are locally minimized, there may exist feasible points in
other regions of the space of decision variables.

 Now we give a simple example to explain the above concepts. Consider the simple quadratic programming problem
\begin{equation}\label{ep-1}
\begin{array}{rl}
 \displaystyle\min & x_1^2+x_2^2\\
{\rm s.t.}& x_1+x_2 -1\leq 0,\\
&-x_1-x_2+2\leq 0.
\end{array}
\end{equation}
It is easy to find that the feasible region is empty.
We consider the minimization problem over the set of points with least constraint violation.
 We can regard the least violation for the constraints as the optimal
value of the following problem
\begin{equation}\label{ep-2}
\begin{array}{rl}
 \displaystyle\min & \displaystyle \frac 12(y_1^2+y_2^2)\\
{\rm s.t.}& x_1+x_2 -1+y_1\leq 0,\\
&-x_1-x_2+2+y_2\leq 0.
\end{array}
\end{equation}
Then the set of points with the least violation is given by
$$
S=\{x: (x,y) \mbox{ solves Problem (\ref{ep-2})}\}.
$$
It is not difficult to obtain
$$
S=\{(x_1,x_2): x_1+x_2-3/2=0\}.
$$
Therefore the minimum point of the objective over the violation set is $(3/4,3/4)$.
It is not difficult to verify that $\varrho(t)=\displaystyle \frac{1}{2}t^2$ for $t\geq 0$ and
$$
\theta (x)=\displaystyle \frac{1}{2} \left\{
[x_1+x_2 -1]_+^2+[-x_1-x_2+2]_+^2
\right\}
$$
is an infeasibility measure in Definition \ref{def:olv} and
$${\rm Argrmin}_x \theta (x)=S.
$$

A popular method for the minimization with least violation constraint is the penalty method. Define  the penalty function as
$$
P_c(x)=x_1^2+x_2^2+c\left\{[x_1+x_2-1]_+^2+[-x_1-x_2+2]_+^2\right\},
$$
where $[t]_+=\max\{0,t\}$ for $t \in \Re$. Obviously $P_c$ is a smooth convex function with
$$
\nabla P_c(x_1,x_2)=\left[
\begin{array}{c}
2x_1+2c\left\{[x_1+x_2-1]_+-[-x_1-x_2+2]_+\right\}\\[3pt]
2x_2+2c\left\{[x_1+x_2-1]_+-[-x_1-x_2+2]_+\right\}
\end{array}
\right].
$$
By solving $\nabla P_c(x)=0$, we obtain
the minimizer of $P_c(x)$ is
\begin{equation}\label{xc}
x^*(c)=\left(\displaystyle \frac{3c}{1+4c},\displaystyle \frac{3c}{1+4c}\right)^T.
\end{equation}
Thus we have
$$
\lim_{c \rightarrow \infty} x^*(c)=\left(\frac 34,\, \frac 34\right)^T;
$$
namely, the limit of the minimizer of $P_c$ approaches the optimal solution. However, for finite $c>0$, the  minimizer of $P_c$ never coincides with the optimal solution. Therefore, it is significant to find a method, different from the penalty method, for solving the minimization optimization problem over the set of points with least violation constraint.

The rest of this paper is organized as follows. In Section 2, we formulate the
minimization problem with the smallest constraint violation for
a general conic optimization problem as a Lipschitz continuous equality constrained optimization problem. For a convex conic problem, we  prove that $\theta(x)$ is convex  and the minimization problem over the set of the least constraint violation can be reformulated as an MPEC problem.
In Section 3, for nonlinear programming with possible inconsistent constraints, we present various types of stationary points for the MPCC problem associated with the minimization problem over the set of the least constraint violation. Especially, an elegant necessary optimality condition, named as L-stationary condition, is established from the classical optimality theory of Lipschitz continuous optimization.
In Section 4, we propose  the smoothing  Fischer-Burmeister function method for  solving the  minimization problem over the set of the least constraint violation for nonlinear programming. It is demonstrated that,
when the positive smoothing parameter approaches to $0$,  any point in the outer limit of
the KKT-point mapping is an L-stationary point of our problem. Some discussions are made in the last section.

\section{Conic Optimization Problem with Least Constraint Voilation}
\setcounter{equation}{0}
In this section, we consider the general mathematical programming problem of the form
\begin{equation}\label{cp}
\begin{array}{ccl}
\mbox{}\hspace*{-6.5cm}({\rm P})\mbox{}\hspace{5cm}\,\ & \displaystyle\min_x &f(x)\\
 & {\rm s.t.}& g(x)\in K,
\end{array}
\end{equation}
where ${\cal X}$  and ${\cal Y}$ are finite dimensional Hilbert  spaces,
$f:{\cal X}\rightarrow \Re$ and $g:{\cal X} \rightarrow {\cal Y}$ are continuously differentiable
mappings, and $K\subset {\cal Y}$
 is a closed convex cone.

For $\varrho(t)=\displaystyle \frac{1}{2}t^2$ for $t\geq 0$, the least violation for the constraint is defined as the optimal
value of the following problem
\begin{equation}\label{ep-cp2}
\begin{array}{cl}
 \displaystyle\min & \displaystyle \frac{1}{2}\|y\|^2\\
{\rm s.t.}& g(x)+y \in K.
\end{array}
\end{equation}
The set of points with the least violation is given by
$$
S=\{x: (x,y) \mbox{ solves Problem (\ref{ep-cp2})}\}.
$$
Our problem is to minimize $f$ over $S$; namely,
\begin{equation}\label{ep-cp2}
\begin{array}{cl}
 \displaystyle\min & f(x)\\[6pt]
{\rm s.t.}& (x,y) \mbox{ solves }\\[6pt]
&\quad  \left\{\begin{array}{cl}
 \displaystyle\min_{w,z} & \displaystyle \frac{1}{2}\|z\|^2\\[4pt]
{\rm s.t.}& g(w)+z \in K.
\end{array}
\right.
\end{array}
\end{equation}
Denote the lower problem of Problem (\ref{ep-cp2}) by ${\rm P_L}$; namely,
\begin{equation}\label{cpL}
\begin{array}{rl}
\mbox{}\hspace*{-6.5cm}({\rm P_L})\mbox{}\hspace{5cm}\,\ \displaystyle \displaystyle\min_{w,z} & \displaystyle \frac{1}{2}\|z\|^2\\[4pt]
{\rm s.t.}& g(w)+z \in K.
\end{array}
\end{equation}
Associated with the above $\varrho (t)$, we have that
\begin{equation}\label{eq:p}
\theta(w)=\displaystyle\min_{z} \left\{ \displaystyle \frac{1}{2}\|z\|^2:
 g(w)+z \in K\right\}.
\end{equation}
It is easy to verify that $\theta (w)$ is an infeasibility measure for Problem (\ref{cp}).
Obviously, we have
$$
\begin{array}{rcl}
\theta(w)& =&\displaystyle\min_{z} \left\{ \displaystyle \frac{1}{2}\|z\|^2:
 g(w)+z \in K\right\}\\[6pt]
 &=&\displaystyle\min_{z'} \left\{ \displaystyle \frac{1}{2}\|z'-g(w)\|^2:
 z' \in K\right\}\\[6pt]
 &=&\displaystyle \frac{1}{2}\|g(w)-\Pi_K(g(w))\|^2.
 \end{array}
$$
Then the optimal value of Problem (${\rm P_L}$)  can be expressed as
$$
{\rm Argmin}\,({\rm P_L})={\rm Argmin}_x \theta(x)={\rm Argmin}_x\displaystyle \frac{1}{2}\|g(x)-\Pi_K(g(x))\|^2={\rm Argmin}_x\displaystyle \frac{1}{2}\|\Pi_{K^{\circ}}(g(x))\|^2,
$$
where $K^{\circ}$ means the polar cone of $K$.
Therefore Problem (\ref{ep-cp2}) can equivalently be expressed as
\begin{equation}\label{simpleE}
\begin{array}{cl}
 \displaystyle\min & \displaystyle f(x)\\
{\rm s.t.}& x \in \mbox{argmin}_w \displaystyle \frac{1}{2}\|\Pi_{K^{\circ}}(g(w))\|^2.
\end{array}
\end{equation}
Now we discuss when the infeasibility measure $\theta$ is a convex function.
\begin{definition}\label{def-convM}
Let $g$ be a continuous mapping. We say that the set-valued mapping ${\cal F}_g: x \rightarrow K-g(x)$ is graph-convex if
$$
\mbox{gph }{\cal F}_g=\{(x,y): g(x)+y \in K\}\subset {\cal X}\times {\cal Y}
$$
is a convex set.
\end{definition}
\begin{proposition}\label{prop:pc}
Let  $g$ be a continuous mapping and ${\cal F}_g: x \rightarrow K-g(x)$ be graph-convex. Then
the function
$$
\theta(x)= \displaystyle \frac{1}{2}\|\Pi_{K^{\circ}}(g(x))\|^2
$$
is convex.
\end{proposition}
{\bf Proof}. From the definition of $\theta(x)$ in (\ref{eq:p}), we have
$$
\theta(x)=\displaystyle\min_{z} \left\{ \displaystyle \frac{1}{2}\|z\|^2:
 g(x)+z \in K\right\}
 = \displaystyle\min_{z} \left\{ \displaystyle \frac{1}{2}\|z\|^2:
(x,z) \in {\rm ghp}\, {\cal F}_g\right\}.
$$
For any $x^i \in {\cal X}$, $i=1,2$, there exist $z^1$ and $z^2$ such that
$$
{\rm Argmin}(x^i)=\displaystyle \frac{1}{2}\|z^i\|^2,\quad (x^i,z^i) \in {\rm ghp}\, {\cal F}_g, i=1,2.
$$
For any $t \in [0,1]$, one has that from the convexity of
${\rm ghp}\, {\cal F}_g$ that
$$
\left((1-t)x^1+tx^2, \, (1-t)z^1+tz^2\right)\in {\rm ghp}\, {\cal F}_g.
$$
Thus
$$
\begin{array}{rcl}
\theta((1-t)x^1+tx^2)&=& \displaystyle\min_{z} \left\{ \displaystyle \frac{1}{2}\|z\|^2:
((1-t)x^1+tx^2,z) \in {\rm ghp}\, {\cal F}_g\right\}\\[6pt]
& \leq & \displaystyle\frac{1}{2}\|(1-t)z^1+tz^2\|^2\\[6pt]
& \leq & (1-t) \displaystyle\frac{1}{2}\|z^1\|^2+t \displaystyle\frac{1}{2}\|z^2\|^2\\[6pt]
&=& (1-t) \theta(x^1)+t \theta(x^2),
\end{array}
$$
which implies that $\theta$ is a convex function.
\hfill $\Box$\\

It follows from the literature that
$\theta(x)$ is differentiable if $K$ is a closed convex cone and  if ${\cal X}=\Re^n$,
\begin{equation}\label{eq:gradent-p}
\nabla \theta(x)={\rm D}g(x)^*\Pi_{K^{\circ}}(g(x)).
\end{equation}
Therefore, for the case when ${\cal X}=\Re^n$, $g:\Re^n \rightarrow Y$ is continuously differentiable, $K \subset Y$ is a closed convex cone, and ${\cal F}_g$ is graph convex, Problem (\ref{simpleE}) is a convex optimization problem and is reduced to
\begin{equation}\label{simpleD}
\begin{array}{cl}
 \displaystyle\min & \displaystyle f(x)\\
{\rm s.t.}& {\rm D}g(x)^*\Pi_{K^{\circ}}(g(x))=0.
\end{array}
\end{equation}
Although Problem (\ref{simpleD}) is a convex optimization problem, we can not handle the constraints easily because they are nonsmooth equalities. We have to transform the constraints to smoothing constraints and then construct numerical algorithms.

If ${\cal F}_g$ defined by Definition \ref{def-convM} is graph-convex, then Problem (${\rm P}_L$) is a convex  optimization problem and thus $(x,y)$ solves (${\rm P}_L$) if and only if there exists $\lambda \in Y^*$ such that
$$
\begin{array}{l}
{\rm D}g(x)^*\lambda=0,\\[4pt]
y+\lambda=0,\\[4pt]
\lambda \in N_K (g(x)+y).
\end{array}
$$
Defining $z=g(x)+y$, the above system can be rewritten as
\begin{equation}\label{eq:cm}
F(x,y,z)=0,\ (y,z) \in \Omega,
\end{equation}
where
\begin{equation}\label{eq:FO}
F(x,y,z)=\left[
\begin{array}{l}
{\rm D}g(x)^*y\\[4pt]
g(x)+y-z
\end{array}
\right] \mbox{ and } \Omega=\{(y,z): K^* \in y \bot z \in K\}.
\end{equation}
Therefore, Problem (\ref{simpleE})  is equivalently expressed as
\begin{equation}\label{simpleC}
\begin{array}{cl}
 \displaystyle\min & \displaystyle f(x)\\
{\rm s.t.}& F(x,y,z)=0,\\
&  (y,z) \in \Omega.
\end{array}
\end{equation}
\section{Convex Nonlinear Programming with Least Constraint Violation}
\setcounter{equation}{0}

For simplicity we consider the following convex nonlinear programming problem
\begin{equation}\label{eq:NLP0}
\begin{array}{ll}
\min & f(x)\\[3pt]
{\rm s.t.} & h(x)=0,\\[3pt]
& g(x)\geq 0,
\end{array}
\end{equation}
where $f:\Re^n\rightarrow \Re$, $h:\Re^n \rightarrow \Re^q$ and $g:\Re^n \rightarrow \Re^p$.
 This is the simple case where $K$ is a closed polyhedral cone:
 \begin{equation}\label{K:simple}
 K=\{0_q\}\times \Re^p_+
 \end{equation}
 and
 $$
 g(x)=\left[
 \begin{array}{c}
 Ax-b\\[4pt]
 c(x)
 \end{array}
 \right],
 $$
 where $A \in \Re^{q\times n}$, $b \in \Re^q$ and $c: \Re^n \rightarrow \Re^p$ with  $c_i(x)$ is concave and twice continuously differentiable for each $i=1,\ldots,p$.
 In this case,
 $$
 F(x,y,z)=\left[
\begin{array}{l}
\displaystyle A^Ty^E+{\cal J}c(x)^Ty^I\\[4pt]
Ax-b+y^E-z^E\\[4pt]
c(x)+y^I-z^I
\end{array}
\right],
  $$
  where $y=(y^E,y^I)$, $z=(z^E,z^I)$ and $y^E=(y_1,\ldots,y_q)^T$, $z^E=(z_1,\ldots,z_q)^T$, $y^I=(y_{q+1},\ldots,y_{q+p})^T$
 and  $z^I=(z_{q+1},\ldots,z_{q+p})^T$. The Jacobian of $F$ at $(x,y,z)$ is of the form
  \begin{equation}\label{JF}
{\cal J} F(x,y,z)=\left[
\begin{array}{ccccc}
\displaystyle \sum_{j=1}^{p}y^I_j\nabla^2 c_j(x) & A^T & {\cal J}c(x)^T & 0 & 0\\[4pt]
A & I_q & 0 & -I_q  & 0\\[4pt]
{\cal J}c(x) & 0 &I_p & 0 & -I_p
\end{array}
\right].
  \end{equation}
 For $K=\{0_q\}\times \Re^p_-$, one has $K^*=\Re^q \times \Re^p_+$ and
 $$
 \Omega=\{(y,x):y^E \in \Re^q, \, z^E=0,\, 0\leq y^I \bot z^I \geq 0\}.
 $$
  We use $\Phi$ to denote the feasible set of Problem
 (\ref{simpleC}); namely,
 $$
 \Phi=\{(x,y,z)\in \Re^n \times \Omega: F(x,y,z)=0\}.
 $$
 Let
$$
\Theta=\{(a, b) \in \Re^p \times \Re^p: 0 \leq a \perp b \geq 0\}.
$$
Then $\Phi$ is simplified  as
\begin{equation}\label{fs}
\Phi=\{(x,y,z): F(x,y,z)=0,\, z^E=0,\, (y^I, z^I)\in \Theta\}
\end{equation}
and Problem (\ref{simpleC}) is simplified as an MPCC problem
\begin{equation}\label{eq:fmpcc}
\min \, f(x) \quad {\rm s.t.}\, (x,y,z)\in \Phi.
\end{equation}

In the following, we derive the tangent cone, the regular normal cone and the normal cone of $\Phi$ at $(x,y, z)\in \Phi$, which are useful in developing S-stationary conditions and M-stationary conditions for Problem (\ref{eq:fmpcc}).

The tangent cone of $\Phi$ at $(x,y, z)$ denoted by $T_{\Phi}(x,y, z)$, the regular normal cone of $\Phi$ at $(x,y, z)$ denoted by $\widehat N_{\Phi}(x,y, z)$ and the normal cone of $\Phi$ at $(x,y,z)$  denoted by $N_{\Phi}(x,y, z)$,
 are defined respectively by
\[
\begin{array}{l}
T_{\Phi}(x,y, z)=\left\{(d_x,d_y,d_z):\left.
 \begin{array}{l}
 \exists t_k \searrow 0, \exists (d^k_x,d^k_y,d^k_z)\rightarrow (d_x,d_y,d_z)\\[3pt]
 \mbox{ satisfying }
(x,y, z)+t_k(d^k_x,d^k_y,d^k_z) \in \Phi
\end{array}
\right.\right\};\\[12pt]
\widehat{N}_{\Phi}(x,y, z)=\left\{(v_x,v_y, v_z):\left.
 \begin{array}{l}
 \langle (v_x,v_y, v_z), (x',y', z') -(x,y, z)\rangle\\[3pt]
\leq {\rm o}(\|(x',y', z') -(x,y, z)\|),\, (x',y', z') \in \Phi
\end{array}\right.
\right\};\\[12pt]
N_{\Phi}(x,y, z)=\left\{(v_x,v_y, v_z):\left.
 \begin{array}{l}
 \exists (x^k,y^k, z^k)\stackrel{\Phi} \rightarrow (x,y, z),
\exists (v^k_x,v^k_y, v^k_z) \rightarrow (v_x,v_y, v_z) \\[3pt]
\mbox{ satisfying }
(v^k_x,v^k_y, v^k_z) \in \widehat N_\Phi (x^k,y^k, z^k)
\end{array}\right.
\right\}.
\end{array}
\]
Let $\omega=\{(\zeta_1,\zeta_2)\in \Re^2_+:\zeta_1\zeta_2=0\}$. For $\Theta$ with  complementarity constraints, we have the following lemma about the variational geometry of $\Theta$ at a point $(\bar a,\bar b) \in \Theta$.
  \begin{lemma}\label{lem:normal-cone}
 For $(\bar a, \bar b)\in
\Theta$, the tangent cone, the regular normal cone and normal cone of $\Theta$ at $(\bar a,\bar b)$ are calculated by
\[T_{\Theta}(\bar a,\bar b)= \bigotimes
_{i=1}^{p} T_{\omega}(\bar a_i,\bar b_i),\,\widehat{N}_{\Omega}(\bar a,\bar b)= \bigotimes
_{i=1}^{p} \widehat{N}_{\omega}(\bar a_i,\bar b_i)
 \mbox{ and }
{N}_{\Theta}(\bar a,\bar b)= \bigotimes _{i=1}^{p}
{N}_{\omega}(\bar a_i,\bar b_i),\] where

\[\bigotimes _{i=1}^{p}
T_{\omega}(\bar a_i,\bar
b_i)=\left\{(u,v)|\,(u_i,v_i)\in T_{\omega}(\bar a_i,\bar b_i),\ i=1,\ldots,p\right\},\]

\[\bigotimes _{i=1}^{p}
\widehat{N}_{\omega}(\bar a_i,\bar
b_i)=\left\{(u,v)|\,(u_i,v_i)\in \widehat{N}_{\omega}(\bar a_i,\bar b_i),\ i=1,\ldots,p\right\},\]
\[\bigotimes _{i=1}^{p}
{N}_{\omega}(\bar a_i,\bar
b_i)=\left\{(u,v)|\,(u_i,v_i)\in {N}_{\omega}(\bar
a_i,\bar b_i),\ i=1,\ldots,p\right\},\]

\[T_{\omega}(\bar a_i,\bar b_i)=\left\{\begin{array}{ll} \Re \times \{0\},& \mbox {if}\  a_i>0,\,b_i=0;\\[6pt]
\{0\}\times \Re, &\mbox {if}\  a_i=0,\,b_i>0;\\[6pt]
\omega,& \mbox {if}\  a_i=0,\,b_i=0,\\
\end{array}\right.\quad\quad\quad\quad\ \quad\quad\quad\quad\quad\quad\quad\]

\[\widehat{N}_{\omega}(\bar a_i,\bar b_i)=\left\{\begin{array}{ll}
\{0\}\times \Re, & \mbox {if}\  a_i>0,\,b_i=0;\\[6pt]
\Re \times \{0\}, &\mbox {if}\  a_i=0,\,b_i>0;\\[6pt]
\Re_-\times\Re_-,& \mbox {if}\  a_i=0,\,b_i=0,\\
\end{array}\right.\quad\quad\quad\quad\ \quad\quad\quad\quad\quad\quad\quad\]

\[{N}_{\omega}(\bar a_i,\bar b_i)=\left\{\begin{array}{ll}\{0\}\times \Re, & \mbox {if}\  a_i>0,\,b_i=0;\\[6pt]
\Re \times \{0\}, & \mbox {if}\  a_i=0,\,b_i>0;\\[6pt]
(\Re \times \{0\})\bigcup(\{0\}\times \Re)\bigcup(\Re_-\times\Re_-), & \mbox {if}\  a_i=0,\,b_i=0.\\
\end{array}\right.\]
 \end{lemma}

For deriving the tangent cone, the regular normal cone and the normal cone of $\Phi$ at  $(x,y, z)\in \Phi$, we need the following assumption.\\[3pt]
{\bf Assumption 1} The Jacobian ${\cal J}F(x,y,z)$ given by (\ref{JF})
is of full row rank.

\begin{proposition}\label{prop-2}
Assume Assumption 1 is satisfied. Then
 \begin{equation}\label{eq-t}
 T_{\Phi}(x,y,z)=\left\{d \in \Re^n \times \Re^{q+p} \times \Re^{q+p}: \begin{array}{l}
  \displaystyle \sum_{j=1}^py^I_j\nabla^2c_j(x)d_x+A^Td_{y^E}+{\cal J}c(x)^Td_{y^I}=0\\[4pt]
  Ad_x+d_{y^E}-  d_{z^E}=0\\[4pt]
    {\cal J}c(x)d_x+d_{y^I}-  d_{z^I}=0\\[4pt]
    d_{z^E}=0\\[4pt]
  (d_{y^I}, d_{z^I}) \in T_\Theta (y^I,z^I)
     \end{array}
   \right
   \},
  \end{equation}
  \begin{equation}\label{eq-n1}
 \widehat N_{\Phi}(x,y,z)=\left
 \{
 \left
 (
 \begin{array}{l}
 \displaystyle \sum_{j=1}^py^I_j\nabla^2c_j(x)\eta_1+
 A^T\eta_2+{\cal J}c(x)^T\eta_3\\
 A\eta_1+\eta_2\\
  {\cal J}c(x)\eta_1+\eta_3+\xi_a\\
 -\eta_2+\eta_4\\
 -\eta_3+\xi_b
 \end{array}
 \right
 ):\begin{array}{l}
 (\eta_1,\eta_2,\eta_3,\eta_4) \in \Re^{n+q+p+q}\\[3pt]
 (\xi_a,\xi_b) \in \widehat N_\Theta (y^I,z^I)
 \end{array}
 \right
 \}
 \end{equation}
 and
 \begin{equation}\label{eq-n2}
  N_{\Phi}(x,y,z)=\left
 \{
 \left
 (
 \begin{array}{l}
 \displaystyle \sum_{j=1}^py^I_j\nabla^2c_j(x)\eta_1+
 A^T\eta_2+{\cal J}c(x)^T\eta_3\\
 A\eta_1+\eta_2\\
  {\cal J}c(x)\eta_1+\eta_3+\xi_a\\
 -\eta_2+\eta_4\\
 -\eta_3+\xi_b
 \end{array}
 \right
 ):\begin{array}{l}
 (\eta_1,\eta_2,\eta_3,\eta_4) \in \Re^{n+q+p+q}\\[3pt]
 (\xi_a,\xi_b) \in  N_\Theta (y^I,z^I)
 \end{array}
 \right
 \}
 \end{equation}
  \end{proposition}
  {\bf Proof.}
 Since   the Jacobian ${\cal J}F(x,y,z)$ given by (\ref{JF})
is of full row rank, we now prove the following equality
 (see 6.7 Exercise of \cite{RockafellarWets1998} for similar result):
 \begin{equation}\label{eq-2}
  T_{\Phi}(x,y,z)=\left\{d: d_{z^E}=0,(d_{y^I},d_{z^I}) \in T_{\Theta}(y^I,z^I):{\cal J}F(x,y,z)d=0\right\}.
 \end{equation}
It is obvious that the set in the left hand-side is contained in the right hand-side and hence we only need to prove the opposite inclusion.
 For any $d=(d_x, d_x, d_z)$ satisfying $d_{z^E}=0,(d_{y^I},d_{z^I}) \in T_{\Theta}(y^I,z^I), {\cal J}F(x,y,z)d=0$, one has that there exist $d^k=(d^x, d^k_y, d^k_z)
 \rightarrow d$ and $t_k \searrow 0$ such that
 $$(x,y^E,z^E,(y^I,z^I))+t_k (d^k_x,d^k_{y^E},d^k_{z^E},(d^k_{y^I},d^k_{z^I})) \in \Re^n\times \Re^q\times \{0_q\}
 \times  \Theta.
 $$
 It follows from Lemma \ref{lem:normal-cone} that $[d_{y^I}]_i [d_{z^I}]_i=0$ for $i=1,\ldots, p$. Let
 \[
 \alpha=\left\{i: y^I_i>0,\, z^I_i=0\right\},\ \beta=\left\{i: y^I_i=z^I_i=0\right\},\ \gamma =\left\{i: y^I_i=0,\, z^I_i>0\right\}
 \]
 and
  \[
 \beta_a=\left\{i\in \beta: [d_{y^I}]_i>0,\, [d_{z^I}]_i=0\right\},\ \beta_b=\left\{i\in \beta: [d_{y^I}]_i=[d_{z^I}]_i=0\right\},\ \beta_c =\left\{i\in \beta: [d_{y^I}]_i=0,\, [d_{z^I}]_i>0\right\}.
 \]
 Let
 \[
 \Gamma_d=\left\{(y^I,z^I)\in \Re^p \times \Re^p:
 \left.
 \begin{array}{l}
 (y^I_{\alpha \cup \beta_a},z^I_{\alpha \cup \beta_a})\in \Re^{|\alpha|+|\beta_a|}_+ \times \{0_{|\alpha|+|\beta_a|}\}\\[6pt]
 (y^I_{\beta_c \cup \gamma},z^I_{\beta_c \cup \gamma})\in \{0_{|\beta_c|+|\gamma|}\}\times \Re^{|\beta_c|+|\gamma|}_+\\[6pt] (y^I_{\beta_b},z^I_{\beta_b})=(0_{|\beta_b|},0_{|\beta_b|})
 \end{array}
 \right.
 \right
 \}.
 \]
 Then $\Gamma_d$ is a convex set and $\Gamma_d \subset \Theta$.
  Since the Jacobian ${\cal J}F(x,y,z)$ given by (\ref{JF})
is of full row rank, it follows from Theorem 2.87 of \cite{BS00} that there exist a neighborhood ${\cal V}$ of $(x,y,z)$ and a positive constant $\kappa$ such that
 \[
 \begin{array}{l}
 {\rm dist}\left((x',[y^E]',[z^E]',([y^I]',[z^I]')),\ [\Re^n\times \Re^q\times \{0_q\}
 \times   \Gamma_d]\cap F^{-1}(0)\right)\\[4pt]
 \quad  \leq \kappa \left \|F(x',y',z'),\, \Pi_{\Re^n\times \Re^q\times \{0_q\}
 \times   \Gamma_d}((x',[y^E]',[z^E]',([y^I]',[z^I]'))\right\|, \,(x',y',z')\in {\cal V}.
 \end{array}
 \]
Noticing that for $(x^k,y^k, z^k)=(x,y, z)+t_k d^k$, $[z^E]^k =0$,
and
\[
\begin{array}{l}
 \left([y^I]^k_{\alpha \cup \beta_a},\, [z^I]^k_{\alpha \cup \beta_a}\right)\in \Re^{|\alpha|+|\beta_a|}_+ \times \{0_{|\alpha|+|\beta_a|}\},\\[6pt]
 \left([y^I]^k_{\beta_c \cup \gamma},\, [z^I]^k_{\beta_c \cup \gamma}\right)\in \{0_{|\beta_c|+|\gamma|}\}\times \Re^{|\beta_c|+|\gamma|}_+,\\[6pt]
 \left([y^I]^k_{\beta_b},\, [z^I]^k_{\beta_b}\right)=\left(t_k[d^k_{y^I}]_{\beta_b},\, t_k[d^k_{z^I}]_{\beta_b}\right )={\rm o}(t_k),
 \end{array}
\]
we have that
\[
\begin{array}{rcl}
 &&\hspace*{-1cm} {\rm dist}\left((x,y, z)+t_k d^k,\, \Phi\right)\\[3pt]
&=&{\rm dist}\left((x,y,z)+t_k d^k,\, [\Re^n\times \Re^q\times \{0_q\}
 \times  \Theta]\cap F^{-1}(0)\right)\\[3pt]
  &=&{\rm dist}\left ((x,y, z)+t_k d^k,\, [\Re^n\times \Re^q\times \{0_q\}
 \times \Gamma_d]\cap F^{-1}(0)\right)\\[3pt]
  &\leq & \kappa \left[\left\|(t_k[d^k_{y^I}]_{\beta_b},\, t_k[d^k_{z^I}]_{\beta_b})\right\|+\left\|F(x^k,y^k, z^k)\right\|\right]\\[5pt]
  &= & \kappa \left\|F(x,y, z)+t_k{\cal J}F(x,y, z)d^k+\int_0^1
 [{\cal J}F((x,y, z)+st_kd^k)-{\cal J}F(x,y, z)]{\rm d}st_kd^k]\right \|\\[3pt]
 &&\quad + \kappa \left\|(t_k[d^k_{y^I}]_{\beta_b},t_k[d^k_{z^I}]_{\beta_b})\right \|\\[3pt]
 &=& {\rm o}(t_k),
 \end{array}
\]
 which implies that $d \in  T_{\Phi}(x,y,z)$. Therefore we obtain equality (\ref{eq-2}).

 Combining with (\ref{JF}) and (\ref{eq-2}), we obtain (\ref{eq-t}).

 Since  ${\cal J}F(x,y,z)$ is of full row rank, formula (\ref{eq-n1}) comes  from the equality
 $$\widehat N_\Phi(x,y,z)
 ={\cal J}F(x,y,z)^*(\Re^n\times \Re^{q+p}\times \Re^{q+p})+\widehat N_{ \Re^n\times \Omega}(x,y,z)
 $$
and
$$
\widehat N_{\Re^n\times \Omega)}(x,y,z))=\{0_n\}\times \{(v_y,v_z):v_{y^E}=0_q, v_{z^E}\in \Re^q,\, (v_{y^I},v_{z^I})\in
\widehat N_\Theta(y^I,z^I)\}.
$$
Formula (\ref{eq-n2}) can be established in the same way. As ${\cal J}F(x,y,z)$ is of full row rank, one has
$$  N_\Phi(x,y,z)
 ={\cal J}F(x,y,z)^*(\Re^n\times \Re^{q+p}\times \Re^{q+p})+ N_{ \Re^n\times \Omega}(x,y,z)
 $$
and
$$
N_{\Re^n\times \Omega)}(x,y,z))=\{0_n\}\times \{(v_y,v_z): v_{y^E}=0_q, v_{z^E}\in \Re^q,\, (v_{y^I},v_{z^I})\in
 N_\Theta(y^I,z^I)\}.
$$
  The proof is completed.  \hfill $\Box$\\

From the above lemma, we can easily develop the necessary optimality conditions for a local minimizer of Problem (\ref{simpleC}). For this purpose, we define
 $$
  \alpha=\left\{i: [y^I]^*_i=0< [z^I]^*_i\right\},\,\,\beta=\left\{i: [y^I]^*_i=0= [z^I]^*_i\right\},\,\,\gamma=\left\{i: [y^I]^*_i>0= [z^I]^*_i\right\}.
  $$

  \begin{theorem}\label{th-S}({\bf S-stationary point})
  Let $(x^*,y^*,z^*)$ be a local minimizer of Problem (\ref{simpleC}). Let Assumption 1 be satisfied at $(x^*,y^*,z^*)$. Then there exist
  $\eta_1^*\in \Re^n$, $[\xi_a]^*_\beta\in \Re^{|\beta|}_-$ and $[\xi_b]^*_{\beta}\in \Re^{|\beta|}_-$ such that
    \begin{equation}\label{s-opta}
  \begin{array}{l}
 \nabla f(x^*)+\left[\displaystyle \sum_{i=1}^p[y^I]^*_j\nabla^2 c_j(x^*)-A^TA-{\cal J}c_{\gamma}(x^*)^T{\cal J}c_{\gamma}(x^*)\right]\eta^*_1+{\cal J}c_{\beta}(x^*)^T[\xi_b]^*_{\beta}=0,\\[8pt]
   {\cal J}c_{\beta}(x^*)\eta^*_1+[\xi_a]^*_{\beta}+[\xi_b]^*_{\beta}=0.
 \end{array}
  \end{equation}
     \end{theorem}
  {\bf Proof}. From \cite{RockafellarWets1998}, we obtain  the inclusion
  \[
  0\in \nabla_{x,y,z}f(x)|_{x=x^*}+\widehat N_{\Phi}(x^*,y^*,z^*),
  \]
  where $\widehat N_{\Phi}(x^*,y^*,z^*)$ is from Lemma \ref{lem:normal-cone}. Namely, there exist  $(\eta_1^*,\eta_2^*,\eta_3^*,\eta_4^*) \in \Re^{n+q+p+q}$ and $([\xi_a]^*,[\xi_b]^*) \in \widehat N_\Theta ([y^I]^*,[z^I]^*)$ such that
  $$
   \left
 (
 \begin{array}{l}
 \nabla f(x^*)+\displaystyle \sum_{j=1}^p[y^I]^*_j\nabla^2c_j(x^*)\eta_1^*+
 A^T\eta_2^*+{\cal J}c(x^*)^T\eta_3^*\\
 A\eta_1^*+\eta_2^*\\
  {\cal J}c(x^*)\eta_1^*+\eta_3^*+[\xi_a]^*\\
 -\eta_2^*+\eta_4^*\\
 -\eta_3^*+[\xi_b]^*
 \end{array}
 \right
 )=0.
 $$
 This set of equalities is simplified as
 \begin{equation}\label{eq:helH}
   \left
 (
 \begin{array}{l}
 \nabla f(x^*)+\displaystyle \sum_{j=1}^p[y^I]^*_j\nabla^2c_j(x^*)\eta_1^*+
 -A^TA\eta_1^*+{\cal J}c(x^*)^T[\xi_b]^*\\
   {\cal J}c(x^*)\eta_1^*+[\xi_a]^*+[\xi_b]^*\\
  \end{array}
 \right
 )=0.
 \end{equation}
 Since  $([\xi_a]^*,\, [\xi_b]^*) \in \widehat N_\Theta ([y^I]^*,\, [z^I]^*)$, one has from Lemma \ref{lem:normal-cone} that ${\cal J}c(x^*)\eta_1^*+[\xi_a]^*+[\xi_b]^*=0$ becomes
 $$
 \begin{array}{l}
 {\cal J}c_{\alpha}(x^*)\eta_1^*+[\xi_a]^*_{\alpha}=0,\ [\xi_b]^*_{\alpha}=0,\\[4pt]
 {\cal J}c_{\gamma}(x^*)\eta_1^*+[\xi_b]^*_{\gamma}=0,\ [\xi_a]^*_{\gamma}=0,\\[4pt]
 {\cal J}c_{\beta}(x^*)\eta_1^*+[\xi_a]^*_{\beta}+[\xi_b]^*_{\beta}=0,\ [\xi_a]^*_{\beta}\leq0,\ [\xi_b]^*_{\beta}\leq 0.
 \end{array}
 $$
 Thus (\ref{eq:helH}) is simplified as (\ref{s-opta}).
   \hfill $\Box$
   \begin{theorem}\label{th-M}({\bf M-stationary point})
  Let $(x^*,y^*,z^*)$ be a local minimizer of Problem (\ref{simpleC}). Let Assumption 1 be satisfied at $(x^*,y^*,z^*)$. Then there exist
  $\eta_1^*\in \Re^n$, $[\xi_a]^*_\beta\in \Re^{|\beta|}$ and $[\xi_b]^*_{\beta}\in \Re^{|\beta|}$ satisfying
  \begin{equation}\label{eq:Msign}
 \Big\{ [\xi_a]^*_i[\xi_b]^*_i=0 \Big\}\mbox{\   or\   } \Big\{ [\xi_a]^*_i\leq 0 \mbox{ and } [\xi_b]^*_i\leq 0 \Big\}\quad \forall i \in \beta,
  \end{equation}
  such that
    \begin{equation}\label{s-optb}
  \begin{array}{l}
 \nabla f(x^*)+\left[\displaystyle \sum_{i=1}^p[y^I]^*_j\nabla^2 c_j(x^*)-A^TA-{\cal J}c_{\gamma}(x^*)^T{\cal J}c_{\gamma}(x^*)\right]\eta^*_1+{\cal J}c_{\beta}(x^*)^T[\xi_b]^*_{\beta}=0,\\[8pt]
   {\cal J}c_{\beta}(x^*)\eta^*_1+[\xi_a]^*_{\beta}+[\xi_b]^*_{\beta}=0.
 \end{array}
  \end{equation}
     \end{theorem}
 Let
\begin{equation}\label{h-mapping}
 H(x,y,z^I)=\left[
\begin{array}{c}
\displaystyle A^Ty^E+{\cal J}c(x)^Ty^I\\[4pt]
Ax-b+y^E\\[4pt]
c(x)+y^I-z^I\\[4pt]
\min\{y^I,z^I\}
\end{array}
\right],
  \end{equation}
 Then  Problem (\ref{simpleC}) is expressed as
\begin{equation}\label{PmC}
\begin{array}{cl}
\displaystyle\min_{x,y,z^I} & f(x)\\
{\rm s.t.}& H(x,y,z^I)=0.
\end{array}
\end{equation}
Noticing that $F$ is a Lipschitz continuous mapping, Problem (\ref{PmC}) is a Lipschitz continuous optimization problem. So we may use the optimality conditions for Lipschitz continuous optimization developed in Clarke (1983). This leads to the so-called C-stationary point. We say that the point $(x^*,y^*,z^*)$ is a C-stationary point if there exist
  $\eta_1^*\in \Re^n$, $[\xi_a]^*_\beta\in \Re^{|\beta|}$ and $[\xi_b]^*_{\beta}\in \Re^{|\beta|}$ satisfying
  \begin{equation}\label{eq:Csign}
[\xi_a]^*_i[\xi_b]^*_i\geq 0 \quad \forall i \in \beta,
  \end{equation}
  such that
    \begin{equation}\label{s-optc}
  \begin{array}{l}
 \nabla f(x^*)+\left[\displaystyle \sum_{i=1}^p[y^I]^*_j\nabla^2 c_j(x^*)-A^TA-{\cal J}c_{\gamma}(x^*)^T{\cal J}c_{\gamma}(x^*)\right]\eta^*_1+{\cal J}c_{\beta}(x^*)^T[\xi_b]^*_{\beta}=0,\\[8pt]
   {\cal J}c_{\beta}(x^*)\eta^*_1+[\xi_a]^*_{\beta}+[\xi_b]^*_{\beta}=0.
 \end{array}
  \end{equation}

It follows from Theorem \ref{th-S} that under Assumption 1,
the point $(x^*,y^*,z^*)$ is a strong stationary point of Problem (\ref{simpleC}). Thus,
 $(x^*,y^*,z^*)$ is an M-stationary point and also a C-stationary point of Problem (\ref{simpleC}).\\

C-stationary conditions are preferred to describe the necessary conditions for MPCC problems. In the following, we will see that for Problem (\ref{PmC}), we may obtain a better  result than conditions in (\ref{s-optc}) by using the optimality conditions for Lipschitz continuous optimization given by
\cite{SchW2007}.

 \begin{proposition}\label{prop-C}({\bf Fritz-John stationary point})
  Let $(x^*,y^*,z^*)$ be a local minimizer of Problem (\ref{PmC}).  Then there exist
  $\eta_0^*\in \Re_+,\eta_1^*\in \Re^n$ and $[v_b]_{\beta}\in \Re^{|\beta|}$ satisfying
  \begin{equation}\label{eq:Msigncc}
[v_b]_i \in [0,1],\ i\in {\beta}
  \end{equation}
  such that
    \begin{equation}\label{s-optcc}
  \begin{array}{l}
 \eta_0^*\nabla f(x^*)+\left[\displaystyle \sum_{i=1}^p[y^I]^*_j\nabla^2 c_j(x^*)-A^TA-{\cal J}c_{\gamma}(x^*)^T{\cal J}c_{\gamma}(x^*)-{\cal J}c_{\beta}(x^*)^T{\rm Diag}([v_b]_{\beta}){\cal J}c_{\beta}(x^*)\right]\eta^*_1=0,
 \end{array}
   \end{equation}
   where ${\rm Diag}(v)={\rm Diag}(v_1,\cdots, v_m)$ for $v \in \Re^m$.
     \end{proposition}
     {\bf Proof}. The generalized Lagrangian of (\ref{PmC}) is
     $$\begin{array}{rl}
     &\hspace{-3mm} L^g(x,y^E,y^I,z^I,\eta_0,\eta_1,\eta_2,\eta_3,\xi)\\[2pt]
     &=\eta_0f(x)+\langle \eta_1, A^Ty^E+{\cal J}c(x)^Ty^I\rangle+
\langle \eta_2, Ax-b+y^E\rangle+\langle \eta_3,
c(x)+y^I-z^I\rangle+ \langle \xi,
\min\{y^I,z^I\}\rangle. \end{array}
     $$
     It follows from the necessary optimality conditions for Lipschitz continuous optimization in \cite{SchW2007} that there exist nonzero vectors $(\eta_0^*,\eta_1^*,\eta_2^*,\eta_3^*,\xi^*)$ with $\eta_0^*\geq 0$ such that
     $$
     0\in \partial_c L\left(x^*,[y^E]^*,[y^I]^*,[z^I]^*,\eta_0^*,\eta_1^*,\eta_2^*,\eta_3^*,\xi^*\right),
     $$
     where $\partial_c$ is the Clarke generalized Jacobian.  Noting that
     \begin{equation}\label{eq:minC}
     \partial_c \min\left\{ [y^I]^*, [z^I]^*\right\}=\left[{\rm Diag}(v_a)  \quad {\rm Diag}(v_b)\right],
     \end{equation}
     where $v_a \in \Re^p$ and $v_b \in \Re^p$ satisfy
     \begin{equation}\label{eq:vab}
     \begin{array}{lll}
     [v_a]_i=1, & [v_b]_i=0, & \mbox{ if } i \in \alpha;\\[5pt]
     [v_a]_i=0, & [v_b]_i=1, & \mbox{ if } i \in \gamma;\\[5pt]
     [v_a]_i=t, & [v_b]_i=1-t, & \mbox{ for some } t\in [0,1] \mbox{ if } i \in \beta.
     \end{array}
     \end{equation}
     Then we get from $0\in \partial_c L\left(x^*,[y^E]^*,[y^I]^*,[z^I]^*,\eta_0^*,\eta_1^*,\eta_2^*,\eta_3^*,\xi^*\right)$ that there exist $v_a \in \Re^p$ and $v_b \in \Re^p$ satisfying (\ref{eq:vab}) such that
     $$
   \left
 (
 \begin{array}{l}
 \eta_0^*\nabla f(x^*)+\displaystyle \sum_{j=1}^p[y^I]^*_j\nabla^2c_j(x^*)\eta_1^*+
 A^T\eta_2^*+{\cal J}c(x^*)^T\eta_3^*\\[6pt]
 A\eta_1^*+\eta_2^*\\[6pt]
  {\cal J}c(x^*)\eta_1^*+\eta_3^*+{\rm Diag}(v_a)\xi^*\\[6pt]
  -\eta_3^*+{\rm Diag}(v_b)\xi^*
 \end{array}
 \right
 )=0.
 $$
 This set of equalities can be simplified as
    \begin{equation}\label{eq:simpx}
   \left
 (
 \begin{array}{l}
 \eta_0^*\nabla f(x^*)+\displaystyle \sum_{j=1}^p[y^I]^*_j\nabla^2c_j(x^*)\eta_1^*-
 A^TA\eta_1^*+{\cal J}c(x^*)^T{\rm Diag}(v_b)\xi^*\\[6pt]
  {\cal J}c(x^*)\eta_1^*+{\rm Diag}(v_b)\xi^*+{\rm Diag}(v_a)\xi^*
  \end{array}
 \right
 )=0.
 \end{equation}
 In view of (\ref{eq:vab}), we have
 $$
 {\rm Diag}(v_b)+{\rm Diag}(v_a)=I_p.
 $$
 So we get from (\ref{eq:simpx}) that $\xi^*=-{\cal J}c(x^*)\eta_1^*$. Substituting this expression back to the first equation
 in (\ref{eq:simpx}), we obtain
 $$
 \eta_0^*\nabla f(x^*)+\displaystyle \sum_{j=1}^p[y^I]^*_j\nabla^2c_j(x^*)\eta_1^*-
 A^TA\eta_1^*-{\cal J}c(x^*)^T{\rm Diag}(v_b){\cal J}c(x^*)=0.
 $$
 Using (\ref{eq:vab}) again, we obtain (\ref{s-optcc}) where $[v_b]_\beta$ satisfies (\ref{eq:Msigncc}).\hfill $\Box$

     From Proposition \ref{prop-C}, we obtain an elegant set of  necessary optimality conditions as follows.
     \begin{theorem}\label{theorem-C}
  Let $(x^*,y^*,z^*)$ be a local minimizer of Problem (\ref{PmC}). Suppose that
  the matrix
  \begin{equation}\label{eq:matrixN}
  \left[\displaystyle \sum_{i=1}^p[y^I]^*_j\nabla^2 c_j(x^*)-A^TA-{\cal J}c_{\gamma}(x^*)^T{\cal J}c_{\gamma}(x^*)\right]
  \end{equation}
  is negatively definite.
   Then there exist
  $\eta_1^*\in \Re^n$ and $[v_b]_{\beta}\in \Re^{|\beta|}$ satisfying
  \begin{equation}\label{eq:Msigncc}
[v_b]_i \in [0,1],\ i\in {\beta}
  \end{equation}
  such that
    \begin{equation}\label{s-optccTa}
  \begin{array}{l}
 \nabla f(x^*)+\left[\displaystyle \sum_{i=1}^p[y^I]^*_j\nabla^2 c_j(x^*)-A^TA-{\cal J}c_{\gamma}(x^*)^T{\cal J}c_{\gamma}(x^*)-{\cal J}c_{\beta}(x^*)^T{\rm Diag}([v_b]_{\beta}){\cal J}c_{\beta}(x^*)\right]\eta^*_1=0.
 \end{array}
   \end{equation}
      \end{theorem}
      \begin{definition}\label{def:Lsp}
      We say that $x^*$ is a L-stationary point for Problem (\ref{PmC}) if there exists $\eta^*_1$ such that
       (\ref{s-optccTa}) is satisfied and (\ref{s-optccTa}) is called L-stationary condition.
      \end{definition}
      Let us introduce the following notation
      \begin{equation}\label{eq:notaS}
      {\cal S}^*=\left\{(x^*,\lambda^*)\in \Re^{2n}:\begin{array}{l}
      \nabla f(x^*)+\bigg[\displaystyle \sum_{i=1}^p[y^I]^*_j\nabla^2 c_j(x^*)-A^TA-{\cal J}c_{\gamma}(x^*)^T{\cal J}c_{\gamma}(x^*) \\[4pt]
      -{\cal J}c_{\beta}(x^*)^T{\rm Diag}([v_b]_{\beta}){\cal J}c_{\beta}(x^*)\bigg]\lambda^*=0
      \end{array}\right\}.
      \end{equation}

      In the next section, we will propose a smoothing function method to generate a sequence of $\{(x^k,\lambda^k)\}$, whose any cluster point is an element of ${\cal S}^*$.
           \begin{remark}\label{th-cCQ}
      The condition that matrix (\ref{eq:matrixN}) is negatively definite is not strict because it  holds if either $A^TA$ is positively definite or $\nabla^2 c_i(x^*)$ is negatively definite and $[y^I]^*_i>0$ for some index $i$.
      \end{remark}
       \begin{remark}\label{comps}
It follows from Theorem \ref{th-S} that under Assumption 1,
the point $(x^*,y^*,z^*)$ is a strong stationary point of Problem (\ref{simpleC}). Thus,
 $(x^*,y^*,z^*)$ is an M-stationary point, also a C-stationary point of Problem (\ref{simpleC}). However, we do not know the relations between Theorem \ref{th-S} and Theorem \ref{theorem-C} because they adopt different constraint qualifications and the results are also different.
 \end{remark}

      \begin{remark}\label{remark-so}
 If there are no inequality constraints, we may check that (\ref{s-optccTa}) or (\ref{s-opta}) is a set of sufficient optimality conditions for $x^*$ being an optimal solution to Problem (\ref{simpleD}) when $f$ is a convex function.
 \end{remark}
\section{The Smoothing  Fischer-Burmeister
Function Method}
\setcounter{equation}{0}
In this section, we only consider the case when the constraints in Problem (\ref{eq:NLP0}) are inconsistent. In this case, the optimization problem with least constraint violation is equivalent to
Problem (\ref{eq:fmpcc}) if $h$ is an affine mapping and $g_i$ is a smooth concave function for $i=1,\ldots, p$. We will present a smoothing function method to solve Problem (\ref{eq:fmpcc}).  Let
\begin{equation}\label{g-mapping}
 G(x,y,z^I)=\left[
\begin{array}{c}
\displaystyle A^Ty^E+{\cal J}c(x)^Ty^I\\[4pt]
Ax-b+y^E\\[4pt]
c(x)+y^I-z^I
\end{array}
\right],
  \end{equation}
Then Problem (\ref{eq:fmpcc}) is an MPCC problem of the following form
\begin{equation}\label{MPCC}
\begin{array}{cl}
\displaystyle\min_{x,y,z^I} & f(x)\\
{\rm s.t.}& G(x,y,z^I)=0,\\[5pt]
& (y^I,z^I)\in \Theta.
\end{array}
\end{equation}
It is well known that, for such a problem, it is not suitable to treat it as a
traditional NLP problem because, as explained in \cite[Example 3.1.1
and Example 3.1.2]{LPRalph96},  even the basic constraint
qualification  (namely, the tangent cone is equal to the linearized cone at an optimal solution) does not hold.
To  overcome  this  difficulty,
various relaxation approaches have been proposed dealing with  the
complementarity constraints.
Facchinei et al. (1999) \cite{FJQi99} and Fukushima and Pang (1999) \cite{FPang99} used
$\psi_{\varepsilon}(a,b)=0$ to approximate the complementarity
relation: $0 \leq a$, $0 \leq  b$, $ab=0$, where
$\psi_{\varepsilon}(a,b)$ is the smoothing  Fischer-Burmeister
function
\begin{equation}\label{eq:sBF}
\psi_{\varepsilon}(a,b)=a+b-\sqrt{a^2+b^2+2\varepsilon^2}.
\end{equation}
Other relaxations of the complementarity relation can be found in for example
Scholtes \cite{Scholtes01}, which uses
$$
a \geq 0,\, b \geq 0,\, ab \leq \varepsilon,
$$
and Lin and Fukushima \cite{LFukushima}, which uses
$$
(a+\varepsilon)(b+\varepsilon) \geq \varepsilon^2 \mbox{ and } ab
\leq \varepsilon^2.
$$
In this section, we shall use
$\psi_{\varepsilon}(a,b)=0$ to approximate the complementarity
relation, where
$\phi_{\varepsilon}(a,b)$ is the smoothing  Fischer-Burmeister
function defined by (\ref{eq:sBF}).

Define
\begin{equation}\label{eq:psim}
\Psi_{\varepsilon}(y^I,z^I)=\left [
\begin{array}{c}
\psi_{\varepsilon}(y^I_1,z^I_1)\\[3pt]
\vdots\\
\psi_{\varepsilon}(y^I_p,z^I_p)
\end{array}
\right ]
\end{equation}
and
 \begin{equation}\label{set-ou}
 \Theta (\varepsilon):=\Big\{(y^I,z^I) \in \Re^p \times \Re^p:
 \Psi_{\varepsilon}(y^I,z^I)=0\Big\}.
 \end{equation}
   Then if $(y^I,z^I) \in \Theta (\varepsilon)$, we have
\[
y^I >0,\ z^I >0 \mbox{ and }y^I_iz^I_i=\varepsilon^2, i=1,\ldots,p.
\]
Obviously, $\psi_0 (a,b)=0$ if and only if $0 \leq a, 0 \leq  b,
ab=0$.  Therefore $\Theta (0)=\Theta$.

For any $(y^I,z^I) \in \Re^{2p}$, we have
\[
 {\cal J}_{y^I,z}\Psi_{\varepsilon}(y^I,z^I)=\left[{\cal J}_{y^I} \Psi_{\varepsilon}(y^I,z^I)\quad {\cal
 J}_{z^I}\Psi_{\varepsilon}(y^I,z^I)\right],
 \]
where
\[
{\cal J}_{y^I} \Psi_{\varepsilon}(y^I,z^I)=
\left
[
\begin{array}{ccc}
1-\displaystyle \frac{[y^I]_1}{\sqrt{[y^I]_1^2+[z^I]_1^2+2\varepsilon^2}} &  &\\[4pt]
& \ddots &\\[4pt]
& & 1-\displaystyle \frac{[y^I]_p}{\sqrt{[y^I]_p^2+[z^I]_p^2+2\varepsilon^2}}
\end{array}
\right ]
\]
and
\[
{\cal J}_{z^I}\Psi_{\varepsilon}(y^I,z^I)= \left [
\begin{array}{ccc}
1-\displaystyle \frac{[z^I]_1}{\sqrt{[y^I]_1^2+[z^I]_1^2+2\varepsilon^2}} &  &\\[4pt]
& \ddots &\\[4pt]
& & 1-\displaystyle \frac{[z^I]_p}{\sqrt{[y^I]_p^2+[z^I]_p^2+2\varepsilon^2}}
\end{array}
\right ].
\]

Let $(y^I, z^I) \in \Theta_\varepsilon$. Then for $i=1,\ldots, p$,
\[
[y^I]_i+[z^I]_i-\sqrt{[y^I]_i^2+[z^I]_i^2+2\varepsilon^2}=0,
\]
we have $[y^I]_i>0$, $[z^I]_i>0$ and $[y^I]_i [z^I]_i=\varepsilon^2$. Thus
$$
\begin{array}{rcl}
1-\displaystyle \frac{[y^I]_i}{\sqrt{[y^I]_i^2+[z^I]_i^2+2\varepsilon^2}}
&=&1-\displaystyle \frac{[y^I]_i}{\sqrt{[y^I]_i^2+[z^I]_i^2+2[y^I]_i [z^I]_i}}\\[6pt]
&=&1-\displaystyle \frac{[y^I]_i}{[y^I]_i+[z^I]_i}\\[6pt]
&=&\displaystyle \frac{[z^I]_i}{[y^I]_i+[z^I]_i},
\end{array}
$$
and in turn we obtain
\begin{equation}\label{z-exp}
1-\displaystyle \frac{[y^I]_i}{\sqrt{[y^I]_i^2+[z^I]_i^2+2\varepsilon^2}}=\displaystyle \frac{[z^I]_i}{[y^I]_i+[z^I]_i},\
1-\displaystyle \frac{[z^I]_i}{\sqrt{[y^I]_i^2+[z^I]_i^2+2\varepsilon^2}}=\displaystyle \frac{[y^I]_i}{[y^I]_i+[z^I]_i}.
\end{equation}

Obviously, for any $\varepsilon >0$, both ${\cal J}_y^I\Psi_{\varepsilon}(y^I,z)$ and
${\cal J}_z\Psi_{\varepsilon}(y^I,z)$ are nonsingular matrices. We can easily
obtain the following conclusion.
 \begin{lemma}\label{lem-LICQ}
 Let $\varepsilon >0$.  Then for any $(y^I,z^I) \in  \Theta (\varepsilon)$,  the linear independence  constraint
 qualification (LICQ) holds and the tangent cone of $\Theta (\varepsilon)$ at
 $(y^I,z^I)$ is
 \begin{equation}\label{eq-tangentcone}
T_{\Theta (\varepsilon)}(y^I,z^I)=\left \{(\triangle y^I,\triangle z^I) \in \Re^{2m}: {\cal J}_{y^I,z^I}\Psi_{\varepsilon}(y^I,z^I)(\triangle y^I,\triangle
z^I)=0\right \},
 \end{equation}
and the normal cone of $\Theta (\varepsilon)$ at
 $(y^I,z^I)$ is
 \begin{equation}\label{eq-normalcone}
N_{\Theta (\varepsilon)}(y^I,z^I)=\widehat N_{\Theta (\varepsilon)}(y^I,z^I)={\cal J}_{y^I,z^I}\Psi_{\varepsilon}(y^I,z^I)^T\Re^p.
 \end{equation}
 \end{lemma}

We use the following problem, denoted by ${\rm P}_\varepsilon$, to approximate Problem (\ref{MPCC}):
\begin{equation}\label{ainlp2}
\begin{array}{cl}
\displaystyle\min_{x,y,z^I} & f(x)\\
{\rm s.t.}& G(x,y, z^I)=0,\\
& (y^I,z^I) \in  \Theta(\varepsilon),
\end{array}
\end{equation}
where $\Theta (\varepsilon)$ is defined by (\ref{set-ou}). Furthermore,
we use $\Phi(\varepsilon)$ to denote the feasible set for Problem (\ref{ainlp2}); namely,
\begin{equation}\label{phi-ep}
\Phi(\varepsilon)=\left\{(x,y^E,y^I, z^I)  \in \Re^n \times \Re^q \times \Theta(\varepsilon): G(x,y,z^I)=0\right\}.
\end{equation}
Define
\begin{equation}\label{Cf}
F_\varepsilon(x,y, z^I)=\left
[
\begin{array}{c}
G(x,y, z^I)\\[4pt]
\Psi_{\varepsilon}(y^I,z^I)
\end{array}
\right
].
\end{equation}
Then $\Phi(\varepsilon)$ is expressed as
\[
\Phi(\varepsilon)=\left\{(x,y^E,y^I, z^I)  \in \Re^n \times \Re^q \times \Re^p \times \Re^p: F_\varepsilon(x,y,z^I)=0\right\}.
\]
By some calculations, we obtain
\begin{equation}\label{eq:JFe}
{\cal J}F_\varepsilon(\vartheta,\mu,y^I,z^I)=\left
[
\begin{array}{cccc}
\displaystyle \sum_{i=1}^p [y^I]_i\nabla^2c_i(x) & A^T & {\cal J}c(x)^T & 0\\[4pt]
A & I & 0 &0\\[4pt]
{\cal J}c(x) & 0 &I &-I\\[4pt]
0 & 0 & {\cal J}_{y^I}\Psi_{\varepsilon}(y^I,z^I) & {\cal J}_{z^I}\Psi_{\varepsilon}(y^I,z^I).
\end{array}
\right].
\end{equation}
Similarly to the proof of Proposition \ref{prop-2}, we can establish the following result.
\begin{proposition}\label{prop-2epsilon}
For $(x,y,z^I)\in \Phi(\varepsilon)$, if
\begin{equation}\label{matrixHx}
\displaystyle \sum_{i=1}^p [y^I]_i\nabla^2c_i(x)-AA^T-{\cal J}c(x)^T{\cal J}_{z^I}\Psi_{\varepsilon}(y^I,z^I){\cal J}c(x)
\end{equation}
is nonsingular, then ${\cal J}F_\varepsilon(x,y,z^I)$ is full of row rank. In this case,
 \begin{equation}\label{eq-tepsilon}
 T_{\Phi(\varepsilon)}(x,y,z^I)
 =\left\{d \in \Re^n \times \Re^{q+p} \times \Re^p: {\cal J}F_\varepsilon(x,y,z^I)d=0\right\}
  \end{equation}
  and
  \begin{equation}\label{eq-nepsilon12}
   N_{\Phi(\varepsilon)}(x,y,z^I)=\widehat N_{\Phi_\varepsilon}(x,y,z^I)
 ={\cal J}F_\varepsilon(\vartheta,\mu,y,z^I)^T\Re^{n+q+p+p}.
  \end{equation}
   \end{proposition}
{\bf Proof}. Let us check that ${\cal J}F_\varepsilon(x,y,z^I)^T$ is of full rank in column. For $\xi_1\in \Re^n$, $\xi_2\in \Re^q$, $\xi_3\in \Re^p$ and  $\xi_4\in \Re^p$, consider
$$
{\cal J}F_\varepsilon(x,y,z^I)^T\left(
\begin{array}{c}
\xi_1\\
\xi_2\\
\xi_3\\
\xi_4
\end{array}
\right
)=0.
$$
It is equivalent to
\begin{equation}\label{eq:checkIn}
\left
[
\begin{array}{l}
\displaystyle \sum_{i=1}^p [y^I]_i\nabla^2c_i(x)\xi_1+A^T\xi_2 +{\cal J}c(x)^T\xi_3\\[4pt]
A\xi_1+\xi_2\\[4pt]
{\cal J}c(x)\xi_1+\xi_3+{\cal J}_{y^I}\Psi_{\varepsilon}(y^I,z^I) \xi_4\\[4pt]
-\xi_3+ {\cal J}_{z^I}\Psi_{\varepsilon}(y^I,z^I)\xi_4.
\end{array}
\right]=0.
\end{equation}
Noting that
$$
{\cal J}_{y^I}\Psi_{\varepsilon}(y^I,z^I)+{\cal J}_{z^I}\Psi_{\varepsilon}(y^I,z^I)
=I_p,
$$
one has from (\ref{eq:checkIn}) that
\begin{equation}\label{checko}
\begin{array}{l}
\xi_2=-A\xi_1\\[4pt]
\xi_4={\cal J}c(x)\xi_1\\[4pt]
\xi_3={\cal J}_{z^I}\Psi_{\varepsilon}(y^I,z^I){\cal J}c(x)\xi_1
\end{array}
\end{equation}
and
\begin{equation}\label{checkd}
\left[\displaystyle \sum_{i=1}^p [y^I]_i\nabla^2c_i(x)-AA^T-{\cal J}c(x)^T{\cal J}_{z^I}\Psi_{\varepsilon}(y^I,z^I){\cal J}c(x)\right]\xi_1=0.
\end{equation}
From the assumption that the matrix of {\ref{matrixHx}) is nonsingular, we obtain from (\ref{checkd}) that $\xi_1=0$ and in turn from  (\ref{checko}) that $\xi_2=0$, $\xi_3=0$ and $\xi_4=0$.  Thus ${\cal J}F_\varepsilon(x,y,z^I)$ is full of row rank and hence
(\ref{eq-tepsilon}) and (\ref{eq-nepsilon12}) follow from Chapter 6 of \cite{RockafellarWets1998}. \hfill $\Box$

\begin{lemma}\label{lem:conv}
For $\Theta(\varepsilon)$ defined by (\ref{set-ou}), we have
\begin{equation}\label{eq:set-connn}
\displaystyle \lim_{\varepsilon \searrow 0} \Theta(\varepsilon)=\Theta (0).
\end{equation}
\end{lemma}
{\bf Proof.} For any $(y^I,z^I) \in \displaystyle \limsup_{\varepsilon \searrow
0} \Theta(\varepsilon)$, there exist $\varepsilon_k \searrow 0$ and $([y^I]^k,[z^I]^k) \in
\Theta (\varepsilon_k)$ such that $([y^I]^k,[z^I]^k)   \rightarrow (y^I,z^I)$. The
inclusion $([y^I]^k,[z^I]^k)  \in \Theta (\varepsilon_k)$ implies that
\[
[y^I]^k+[z^I]^k -\sqrt{([y^I]^k)^2+([z^I]^k)^2+2\varepsilon_k^2}=0.
\]
Then, letting $k  \rightarrow \infty$, we have
\[
y^I+z^I-\sqrt{[y^I]^2+[z^I]^2}=0;
\]
namely, $\psi_0(y^I,z^I)=0$ and $(y^I,z^I) \in \Theta(0)$. Therefore we have
\[
\displaystyle \limsup_{\varepsilon \searrow 0} \Theta(\varepsilon) \subset \Theta
(0).
\]
For any $(y^I,z^I) \in \Theta(0)$, let
\[
I_+=\left\{i: [y^I]_i >0\right\},\ J_+=\left\{i:[z^I]_i >0\right\},\ I_0=\{1,\ldots,m\}\setminus
\Big(I_+\cup J_+\Big).
\]
For any $\varepsilon >0$ defined $(y^I(\varepsilon),z(\varepsilon))$ by
\begin{equation}\label{eq:aux-help}
([y^I]_i(\varepsilon),[z^I]_i(\varepsilon))=\left \{
\begin{array}{ll}
([y^I]_i, \varepsilon^2/[y^I]_i), & \mbox{ if } i \in I_+;\\[4pt]
( \varepsilon^2/[z^I]_i,[z^I]_i), & \mbox{ if } i \in J_+;\\[5pt]
(\varepsilon,\varepsilon) & \mbox{ if } i \in I_0,
\end{array}
\right.
\end{equation}
Then $\psi_{\varepsilon}([y^I]_i(\varepsilon),[z^I]_i(\varepsilon))=0$ for $i=1,\ldots, m$. Thus
$\Psi_{\varepsilon}(y^I(\varepsilon),z^I(\varepsilon))=0$ or, equivalently, $(y^I(\varepsilon),z^I(\varepsilon) \in
\Theta (\varepsilon)$. Obviously, $(y^I(\varepsilon),z^I(\varepsilon) \rightarrow (y^I,z^I)$. This
implies that
\[
\displaystyle \liminf_{\varepsilon \searrow 0} \Theta(\varepsilon) \supset \Theta
(0).
\]
Therefore $\Theta(\varepsilon) \rightarrow \Theta(0)$ as $\varepsilon \searrow 0$.
\hfill $\Box$\\
\begin{corollary}\label{coro-phi-cov}
Let $\Phi(\varepsilon)$ be defined by (\ref{phi-ep}). Then
$$
\Phi (\varepsilon) \rightarrow\Phi\mbox{ as }\varepsilon \searrow 0.
$$
\end{corollary}
{\bf Proof}. The result can be obtained by noting that $\Phi (\varepsilon)$ and $\Phi$ can be expressed as
\[
\Phi(\varepsilon)=\{(x,y^E,y^I, z^I)  \in \Re^n \times \Re^q \times \Re^p\times \Re^p: G(x,y^E,y^I,z^I)=0\}\cap \Re^n \times \Re^q \times   \Theta(\varepsilon)
\]
and
\[
\Phi=\{(x,y^E,y^I, z^I)  \in \Re^n \times \Re^q \times \Re^p\times \Re^p: G(x,y^E,y^I,z^I)=0\}\cap \Re^n \times \Re^q \times   \Theta,
\]
respectively. \hfill $\Box$\\

Now denote the optimal value and the (global) solution set of Problem ${\rm P}_\varepsilon$ by
$\kappa (\varepsilon)$ and $S(\varepsilon)$, respectively; namely,
$$
\begin{array}{l}
 \kappa (\varepsilon):=\inf \{f(x) \,|\,(x,y,z^I) \in \Phi (\varepsilon)\},\\
 \\
 S(\varepsilon):=\mbox{Argmin}\{\{f(x) \,|\,(x,y,z^I) \in \Phi (\varepsilon)\}\}.
 \end{array}
 $$
\begin{theorem}\label{th3.2}
Let $f$ be level-bounded; namely, the level set of $f$ is bounded. Let  ${\rm P}_\varepsilon$ is defined by (\ref{ainlp2}), and $\kappa (\varepsilon)$ and $S(\varepsilon)$ be its optimal value and solution set, respectively.   Then the function $\kappa
(\varepsilon)$ is continuous at $0$ with respect to $\Re_+$ and the
set-valued mapping $S(\varepsilon)$ is outer semi-continuous at $0$ with
respect to $\Re_+$.
\end{theorem}
{\bf Proof.} As $f$ is level-bounded, we have $\kappa (\varepsilon)$ is finite
and $S(\varepsilon) \ne \emptyset$ for any $\varepsilon \geq 0$.
Let
\[
\widehat{f}_{\varepsilon}(x,y, z^I)=f(x)+\delta_{\Phi (\varepsilon)}(x,y, z^I),
\]
where $\delta_{\Phi(\varepsilon)}$ is the indicator function of $\Phi
(\varepsilon)$. From Lemma \ref{lem:conv}, $\Phi (\varepsilon)\rightarrow
\Theta(0)$ as $\varepsilon \searrow 0$, $\widehat{f}_{\varepsilon}$ epi-converges to
$\widehat{f}_{0}$.  The level-boundedness of  $\widehat{f}_{\varepsilon}$ is
easily verified for $\varepsilon \geq 0$. Therefore, we have from Theorem
7.41 of Rockafellar and Wets (1998) that the function $\kappa (\varepsilon)$
is continuous at $0$ with respect to $\Re_+$ and the set-valued
mapping $S(\varepsilon)$ is outer semi-continuous at $0$ with respect to
$\Re_+$. The proof is completed. \hfill $\Box$

If  $(x,y, z^I)\in \Phi(\varepsilon)$ is a local minimizer of ${\rm P}_\varepsilon$  and ${\cal J}F_{\varepsilon}(x,y, z^I)$ is of full row rank, then there exists a vector
$\xi\in \Re^{n+q+2p}$
such that
$$
\nabla_{x,y, z^I} f(x)+{\cal J}F_{\varepsilon}(x,y, z^I)^T\xi=0,
$$
which  is reduced to
$$
\nabla_x f(x)+\left[\displaystyle \sum_{i=1}^p [y^I]_i\nabla^2c_i(x)-AA^T-{\cal J}c(x)^T{\cal J}_{z^I}\Psi_{\varepsilon}(y^I,z^I){\cal J}c(x)\right]\xi_1=0.
$$
This leads to the following definition.
\begin{definition}\label{def:epsilon-stat}
We say $(x,y, z^I)\in \Phi(\varepsilon)$ is a stationary point of ${\rm P}_\varepsilon$ if there exists a vector
$\lambda \in \Re^{n}$
\begin{equation}\label{epsilon-st}
\nabla_x f(x)+\left[\displaystyle \sum_{i=1}^p [y^I]_i\nabla^2c_i(x)-AA^T-{\cal J}c(x)^T{\cal J}_{z^I}\Psi_{\varepsilon}(y^I,z^I){\cal J}c(x)\right]\lambda=0.
\end{equation}

\end{definition}
The following theorem is about the convergence of the stationary points for ${\rm P}_\varepsilon$, which  shows that a cluster point of stationary points for ${\rm P}_\varepsilon$ is related to the condition  (\ref{s-optccT}) when $\varepsilon \searrow 0$.
\begin{theorem}\label{convergence-sm}
Let  $(x (\varepsilon),y(\varepsilon), z^I(\varepsilon))\in \Re^{n+q+2p}$ be a stationary point for ${\rm P}_\varepsilon$ for $\varepsilon >0$, with multiplier $\lambda (\varepsilon)\in \Re^{n}$. Then for any
 \[
 \left(x^*,y^*,[z^I]^*,\lambda^*\right) \in \displaystyle \limsup_{\varepsilon \searrow 0}\left \{(x (\varepsilon),y(\varepsilon), z^I(\varepsilon),\lambda(\varepsilon))\right\},
 \]
 one has that $(x^*,\lambda^*) \in {\cal S}^*$, where ${\cal S}^*$ is defined by (\ref{eq:notaS}).
\end{theorem}
{\bf Proof}. Let $(x^*,y^*, [z^I]^*,\lambda^*) \in
 \displaystyle \limsup_{\varepsilon \searrow 0} \left\{(x (\varepsilon),y(\varepsilon), z^I(\varepsilon),\lambda(\varepsilon))\right\}$.  Then there exists a sequence $\varepsilon_k \searrow 0$  and $\left(x^k,y^k,[z^I]^k,\lambda^k\right)=\left(x(\varepsilon_k),y(\varepsilon_k),[z^I](\varepsilon_k),\lambda(\varepsilon_k)\right)$ such that $\left(x^k,y^k,[z^I]^k,\lambda^k\right)\rightarrow \left(x^*,y^*, [z^I]^*,\lambda^*\right)$ with
 \begin{equation}\label{eq:h11}
 \nabla_x f(x^k)+\left[\displaystyle \sum_{i=1}^p [y^I]^k_i\nabla^2c_i(x^k)-AA^T-{\cal J}c(x^k)^T{\cal J}_{z^I}\Psi_{\varepsilon}([y^I]^k,[z^I]^k){\cal J}c(x^k)\right]\lambda^k=0.
 \end{equation}
  It follows from Lemma \ref{lem:conv} that $([y^I]^*,[z^I]^*)\in \Theta$. Define
 $$
  \alpha=\left\{i: [y^I]^*_i=0< [z^I]^*_i\right\},\,\,\beta=\left\{i: [y^I]^*_i=0=[z^I]^*_i\right\},\,\,\gamma=\left\{i: [y^I]^*_i>0=[z^I]^*_i\right\}.
  $$
 Noting that
 $$
 {\cal J}_{z^I}\Psi_{\varepsilon}([y^I]^k,[z^I]^k) =\left
 [
 \begin{array}{ccc}
 \displaystyle \frac{[y^I]^k_1}{[z^I]^k_1+[y^I]^k_1} & &\\[4pt]
 & \ddots &\\[4pt]
 & &  \displaystyle \frac{[y^I]^k_p}{[z^I]^k_p+[y^I]^k_p}
 \end{array}
 \right
 ],
 $$
 we have
 $$
 \displaystyle \frac{[y^I]^k_i}{[z^I]^k_i+[y^I]^k_i}\rightarrow \left \{
 \begin{array}{ll}
 0, & i \in \alpha;\\[4pt]
 1, & i \in \gamma.
 \end{array}
 \right.
 $$
 For $i \in \beta$, since $\displaystyle \frac{[y^I]^k_i}{[z^I]^k_i+[y^I]^k_i}\in (0,1)$, it has an cluster point $\eta_i \in [0,1]$. Thus there exists $\{k_m:m \in \textbf{N}\}$ such that
  \begin{equation}\label{eq:clusterm}
  \displaystyle \frac{[y^I]^{k_m}_i}{[z^I]^{k_m}_i+[y^I]^{k_m}_i}\rightarrow \left \{
 \begin{array}{ll}
 0, & i \in \alpha;\\[4pt]
 \eta_i, & i \in \beta;\\[4pt]
 1, & i \in \gamma.
 \end{array}
 \right.
 \end{equation}
 Taking the limit for $k=k_m$, $m \rightarrow \infty$ in (\ref{eq:h11}), we obtain
 \begin{equation}\label{s-optccT}
  \begin{array}{l}
 \nabla f(x^*)+\left[\displaystyle \sum_{i=1}^p[y^I]^*_j\nabla^2 c_j(x^*)-A^TA-{\cal J}c_{\gamma}(x^*)^T{\cal J}c_{\gamma}(x^*)-{\cal J}c_{\beta}(x^*)^T{\rm Diag}(\eta_{\beta}){\cal J}c_{\beta}(x^*)\right]\lambda^*=0
 \end{array}
    \end{equation}
with $\eta_i \in [0,1]$ for $i \in \beta$. Thus $(x^*,y^*, [z^I]^*,\lambda^*)$ satisfied (\ref{s-optccT}) and $(x^*,\lambda^*)\in {\cal S}^*$.
 The proof is completed.
 \hfill $\Box$
\begin{theorem}\label{th:convergence}
Let  $(x (\varepsilon),y(\varepsilon), z^I(\varepsilon))\in \Re^{n+q+2p}$ be a local minimizer of ${\rm P}_\varepsilon$ for $\varepsilon >0$. Let
\[
 \left(x^*,y^*,[z^I]^*\right) \in \displaystyle \limsup_{\varepsilon \searrow 0}\left \{\left(x (\varepsilon),y(\varepsilon), z^I(\varepsilon)\right)\right\}.
 \]
 If the matrix
  \begin{equation}\label{eq:matrixNe}
  \left[\displaystyle \sum_{i=1}^p[y^I]^*_j\nabla^2 c_j(x^*)-A^TA-{\cal J}c_{\gamma}(x^*)^T{\cal J}c_{\gamma}(x^*)\right]
  \end{equation}
  is negatively definite, then there exists a vector $\lambda^* \in \Re^n$ such that $(x^*,\lambda^*) \in {\cal S}^*$.
\end{theorem}
{\bf Proof}.
For $\left(x^*,y^*, [z^I]^*\right) \in
 \displaystyle \limsup_{\varepsilon \searrow 0} \left\{(x (\varepsilon),y(\varepsilon), z^I(\varepsilon)))\right\}$, there exists a sequence $\varepsilon_k \searrow 0$  and $\left(x^k,y^k,[z^I]^k\right)=\left(x(\varepsilon_k),y(\varepsilon_k),[z^I](\varepsilon_k)\right)$ such that $\left(x^k,y^k,[z^I]^k\right)\rightarrow \left(x^*,y^*, [z^I]^*\right)$. Since the matrix in (\ref{eq:matrixNe}) is negatively definite, the matrix
 \begin{equation}\label{eq:h11em}
 \left[\displaystyle \sum_{i=1}^p [y^I]^k_i\nabla^2c_i(x^k)-AA^T-{\cal J}c(x^k)^T{\cal J}_{z^I}\Psi_{\varepsilon}([y^I]^k,[z^I]^k){\cal J}c(x^k)\right]
 \end{equation}
is negatively definite for $k$ large enough. Then there exists a unique vector $\lambda^k \in \Re^n$ such that
$$
\nabla f(x^k)+\left[\displaystyle \sum_{i=1}^p [y^I]^k_i\nabla^2c_i(x^k)-AA^T-{\cal J}c(x^k)^T{\cal J}_{z^I}\Psi_{\varepsilon}([y^I]^k,[z^I]^k){\cal J}c(x^k)\right]\lambda^k=0.
$$
Then
$$
\lambda^k=-\left[\displaystyle \sum_{i=1}^p [y^I]^k_i\nabla^2c_i(x^k)-AA^T-{\cal J}c(x^k)^T{\cal J}_{z^I}\Psi_{\varepsilon}([y^I]^k,[z^I]^k){\cal J}c(x^k)\right]^{-1}\nabla f(x^k)
$$
and $\lambda^k$ has a cluster point $\lambda^*$ such that there exists
   $[v_b]_{\beta}\in \Re^{|\beta|}$ satisfying
 $$
[v_b]_i \in [0,1],\ i\in {\beta}
 $$
  and
  $$
  \begin{array}{l}
 \nabla f(x^*)+\left[\displaystyle \sum_{i=1}^p[y^I]^*_j\nabla^2 c_j(x^*)-A^TA-{\cal J}c_{\gamma}(x^*)^T{\cal J}c_{\gamma}(x^*)-{\cal J}c_{\beta}(x^*)^T{\rm Diag}([v_b]_{\beta}){\cal J}c_{\beta}(x^*)\right]\lambda^*=0.
 \end{array}
   $$
The proof is completed. \hfill $\Box$

From the above theorem, we see that the smoothing Fischer-Burmeister function method works for the convex nonlinear programming with inconsistent constraints. Specifically, when the positive smoothing parameter of the method approaches to zero, any point in the outer limit of the KKT-point mapping is an L-stationary point of the equivalent MPCC problem.

\section{Discussions}
This paper established the optimization model with least constraint violation to model constrained optimization problems with possible inconsistent constraints. If the constraints are consistent, the model is reduced to the original problem. When the constraints in a conic optimization problem are possible inconsistent, the model is reformulated as an MPEC problem. For the nonlinear programming problem with possible inconsistent convex constraints,  several stationary points for the equivalent MPCC problem were given. Importantly,  the  so-called L-stationary point is proposed, from the optimality theory for Lipschitz continuous optimization. The  smoothing  Fischer-Burmeister
function method is constructed to solve the equivalent MPCC problem and any accumulation point of the sequence generated by the smoothing function method is an L-stationary point.

There are many  topics  left to investigate for optimization with least constraint violation. When we do not know whether the feasible region is nonempty or not, the optimization problem with least violation is always feasible, this is its advantage. However, if the constrained problem is feasible, then the model involves the infeasibility measure $\theta (x)$, which is usually smooth but not twice differentiable even functions of the original problem are all twice differentiable, this brings computational difficulties. The smoothing Fischer-Burmeister
function method only deals with the case when the constraints are inconsistent for convex nonlinear programming, it has nothing to do with the original problem when it is feasible. Is it possible for us to propose a unified algorithm, which can solve the optimization problem with least violation no matter when the original problem is either infeasible or feasible? Another question is as follows. The smoothing function algorithm can only cope with the nonlinear programming problem with least constraint violation. Can we construct algorithms to deal with other conic optimization problems, for example nonlinear semidefinite optimization problem?

\medskip \noindent
{\bf Acknowledgments.} The authors are very grateful to Professors Ya-xiang Yuan, Xinwei Liu and Zhongwen Chen for their useful discussions and comments.


\begin{thebibliography}{99}
\bibitem{BS00}
 J. F. Bonnans and A. Shapiro, {Perturbation Analysis of Optimization Problems},
New York, Springer, 2000.

\bibitem{Byrd2010}
R. H. Byrd, F. E. Curtis and J. Nocedal,
{\sl Infeasibility Detection and SQP Methods for Nonlinear Optimization},
SIAM J. Optim., 20:5 (2010), 2281-2299.
\bibitem{Burke2014}
J. V. Burke, F. E. Curtis and H. Wang, A sequential quadratic optimization algorithm with
rapid infeasibility detection, SIAM J. Optim., 24 (2014), 839-872.
\bibitem{Clarke83}
F. H. Clarke,  {Optimization and Nonsmooth Analysis},
John Wiley and Sons, New York, 1983.
\bibitem{DLSun2020}
 Y. H. Dai, X. W. Liu and J. Sun,
{\sl A primal-dual interior-point method capable of rapidly detecting infeasibility for nonlinear programs}, Journal of Industrial and Management Optimization, 16:2 (2020), 1009-1035.
\bibitem{FJQi99}
F. Facchinei, H. Jiang  and L. Qi, {\sl  A smoothing
method for mathematical programs with equilibrium constraints},
Math. Prog., 85 (1999), 107-134.

\bibitem{FPang99}
 M. Fukushima and J. S. Pang, {\sl Convergence of a smoothing continuation
method for mathematical problems with complementarity constraints,
Ill-posed Variational Problems and Regularization Techniques},
Lecture Notes in Economics and Mathematical Systems, Vol. 477,
Th$\acute{e}$ra, M. and Tichatschke, R.(eds.), Springer-Verlag,
Berlin/Heidelberg, 1999, 105-116.

\bibitem{LFukushima}
G. H. Lin  and M. Fukushima, {\sl A modified relaxation
scheme for mathematical programs with complementarity constraints},
Ann. Oper. Res., 133 (2005), 63-84.

\bibitem{LPRalph96}
Z. Q. Luo, J. S. Pang and D. Ralph, {Mathematical
Programs with Equilibrium Constraints}, Cambridge University Press, 1996.
\bibitem{NWright99}
J. Nocedal and S. J. Wright, Numerical Optimization, Springer Press, 1999.
\bibitem{Rob80} S. M. Robinson, {\sl Strongly regular generalized equations}, Mathematics of Operations Research 5 (1980), 43-62.
\bibitem{Rob81}
 S. M. Robinson, {\sl Some continuity properties of polyhedral multifunctions}, Mathematical Programming Study, 14 (1981), 206-214.

\bibitem{Rob82} S. M. Robinson, {\sl Generalized Equations and Their Solutions, Part II: Applications to Nonlinear Programming}, Mathematical Programming Study 19 (1982), 200-221.
\bibitem{RockafellarWets1998}  R. T. Rockafellar and R. J. -B. Wets, {Variational Analysis},  Berlin-Heidelberg:  Springer-Verlag, 1998.
\bibitem{Scholtes01}
 S. Scholtes, {\sl  Convergence properties of a regularization scheme for
mathematical programs with complementarity constraints}, SIAM J.
Optim., 11 (2001), 918-936.

\bibitem{SchW2007}
 W. Schirotzek, Nonsmooth Analysis, Berlin Heidelberg: Springer-Verlag, 2007. 
\end{thebibliography}
\end{document}